# New RBF collocation schemes and their applications

W. CHEN

Department of Informatics, University of Oslo, Norway

(Sept. 9, 2001)

**1. Motivations**

**2. Symmetric boundary knot method**

**3. Boundary particle method**

**4. Numerical validations of the BKM and BPM**

5. **Modified Kansa method and its numerical validations**

6. **Remarks: merits and demerits**

# 1. Motivations

The radial basis function (RBF) method is truly meshfree and independent of dimensionality and geometry complicity and has inherent multiscale capability. Among the existing RBF-based schemes for PDE's are

1. Domain-type schemes: Kansa's method (unsymmetric) and Fasshauer's Hermite method (symmetric). Both lose significant accuracy nearby boundary.

2. Boundary-type schemes: method of fundamental solution (MFS), also known as regular BEM. The method is unsymmetric and requires controversial fictitious boundary outside physical domain due to singularity of fundamental solution, which causes instability for irregular geometry.

The purposes of this study are to

1. develop a symmetric boundary knot method (BKM) which bases on the Hermite interpolation with nonsingular general solution and uses the dual reciprocity principle (DRM) to evaluate particular solution;

2. introduce a truly boundary-only boundary particle method (BPM) which applies the multiple reciprocity principle (MRM);

3. present a domain-type modified Kansa method (MKM) by combining symmetric Hermite interpolation and the DRM to improve the solution accuracy close to boundary.

# 2. Symmetric boundary knot method and boundary particle method

The BKM can be viewed as a two-step scheme, approximation of particular solution and the evaluation of homogeneous solution. Let us consider the differential equation

$$\Re\{u\} = f(x), \qquad x \in \Omega$$

with boundary conditions

$$u(x) = R(x), \qquad x \subset S_u,$$

$$\frac{\partial u(x)}{\partial n} = N(x), \qquad x \subset S_T.$$

The solution of the above equation can be split as

$$u = u_h + u_p.$$

The particular solution $u_p$ satisfies the governing equation but not necessarily boundary conditions, while the homogeneous solution $u_h$ must hold both, namely,

$$\Re\{u_h\} = 0,$$

$$u_h(x) = R(x) - u_p,$$

$$\frac{\partial u_h(x)}{\partial n} = N(x) - \frac{\partial u_p(x)}{\partial n}.$$

Like the MFS and BEM, the DRM and RBF are employed to evaluate the particular solution. The inhomogeneous term is approximated by

$$f(x) \cong \sum_{j=1}^{N+L} \alpha_j \phi(r_j).$$

We have

$$u_p = \sum_{j=1}^{N+L} \alpha_j \varphi\left(\|x - x_j\|\right),$$

where the RBF $\phi$ is related to the RBF $\varphi$ through operator $\Re$.

The distinctions of the BKM are to use the nonsingular general solution, namely,

$$u_h(x) = \sum_{k=1}^{L} \lambda_k u_0^\#(\|x - x_k\|)$$

Unlike the MFS, the BKM places all nodes only on physical boundary. However, the naïve use of the above representation leads to an unsymmetric scheme. Instead, we use

$$u_h(x) = \sum_{s=1}^{L_d} \lambda_s u_0^\#(r_s) - \sum_{s=L_d+1}^{L_d+L_N} \lambda_s \frac{\partial u_0^\#(r_s)}{\partial n},$$

Substituting the above RBF representation into boundary equations produces

$$\sum_{s=1}^{L_d} a_s u_0^\#(r_{is}) - \sum_{s=L_d+1}^{L_d+L_N} a_s \frac{\partial u_0^\#(r_{is})}{\partial n} = R(x_i) - u_p(x_i),$$

$$\sum_{s=1}^{L_d} a_s \frac{\partial u_0^\#(r_{js})}{\partial n} - \sum_{s=L_d+1}^{L_d+L_N} a_s \frac{\partial^2 u_0^\#(r_{js})}{\partial n^2} = N(x_j) - \frac{\partial u_p(x_j)}{\partial n},$$

$$\sum_{s=1}^{L_d} a_s u_0^\#(r_{ls}) - \sum_{s=L_d+1}^{L_d+L_N} a_s \frac{\partial u_0^\#(r_{ls})}{\partial n} = u_l - u_p(x_l).$$

# 3. Boundary particle method

According to the multiple reciprocity theorem, the particular solution can be approximated by higher-order homogeneous solution

$$u = u_h^0 + u_p^0 = u_h^0 + \sum_{m=1}^{\infty} u_h^m.$$

Through an incremental differentiation via operator $\Re\{\}$, we have:

$$\sum_{s=1}^{L_d} \beta_s u_0^\#(r_i) - \sum_{s=L_d+1}^{L_d+L_N} \beta_s \frac{\partial u_0^\#(r_i)}{\partial n} = R(x_i) - u_p^0(x_i)$$

$$\sum_{s=1}^{L_d} \beta_s \frac{\partial u_0^\#(r_j)}{\partial n} - \sum_{s=L_d+1}^{L_d+L_N} \beta_s \frac{\partial^2 u_0^\#(r_j)}{\partial n^2} = N(x_j) - \frac{\partial u_p^0(x_j)}{\partial n}$$

$$\sum_{s=1}^{L_d} \beta_s \Re^{n-1}\{u_n^\#(r_{is})\} - \sum_{s=L_d+1}^{L_d+L_N} \beta_s \frac{\partial \Re^{n-1}\{u_n^\#(r_{is})\}}{\partial n} = \Re^{n-2}\{f(x_i)\} - \Re^{n-1}\{u_p^n(x_i)\}$$

$$\sum_{s=1}^{L_d} \beta_s \frac{\partial \Re^{n-1}\{u_n^\#(r_{js})\}}{\partial n} - \sum_{s=L_d+1}^{L_d+L_N} \beta_s \frac{\partial^2 \Re^{n-1}\{u_n^\#(r_{js})\}}{\partial n^2} = \frac{\partial \left( \Re^{n-2}\{f(x_j)\} - \Re^{n-1}\{u_p^n(x_j)\} \right)}{\partial n}$$

$$n=1,2,\ldots,$$

where $\Re^n\{\}$ denotes the *n*-th order operator $\Re\{\}$.

The successive process is truncated at some order $M$. The practical solution procedure is a reversal recursive process:

$$\beta_k^M \to \beta_k^{M-1} \to \cdots \to \beta_k^0.$$

It is noted that due to

$$\Re^{n-1}\{u_h^n(r_k)\} = u_h^0(r_k),$$

the coefficient matrices of all successive equation are the same, i.e.

$$Q\beta^n = b^n, \quad n=M, M-1,\ldots,1,0.$$

Thus, the LU decomposition algorithm is suitable. Finally, we have

$$u(x_i) = \sum_{n=0}^{M}\sum_{k=1}^{L} \beta_k^n u_n^\#(r_{ik}).$$

# 4. Numerical validations for the BKM and BPM

## 4.1. Helmholtz problem

$$\nabla^2 u + \gamma^2 u = f(x)$$

with Dirichlet and Neumann boundary conditions. The exact solutions are

$$u = x^2 \sin(dx)\cos(dy)$$

for 2D inhomogeneous Helmholtz problem ($\gamma = d\sqrt{2}$) and

$$u = \cos(dx)\sin(dy)\sin(dz)$$

for 3D homogeneous Helmholtz ($\gamma = d\sqrt{3}$).

## 4.2. Steady convection-diffusion problem

$$D\nabla^2 u - v \bullet \nabla u - \kappa u = g(x)$$

with Dirichlet and Neumann boundary conditions. The exact solutions are

$$u = x^2 e^{-\eta(x+y)}$$

for 2D inhomogeneous problem, where $D=1$, $v_x=v_y=-\sigma$, $\kappa=3\sigma^2/2$, $\eta = \left(\sigma + \sqrt{\sigma^2 + 2\kappa}\right)/2$, and

$$u = e^{-\sigma(x+y+z)}$$

for 3D homogeneous problem, where $D=1$, $v_x=v_y=v_z=-\sigma$, $\kappa=7\sigma^2/12$.

The $L_2$ norms of relative errors are calculated at 460 nodes for 2D cases and 1012 knots for 3D cases.

Table 1. $L_2$ norm of relative errors for 2D inhomogeneous Helmholtz problems by the BKM and BPM

|  | BKM (41+15) | BKM (49+15) | BPM (49) | BPM (65) |
|---|---|---|---|---|
| $\gamma = \sqrt{2}$ | 1.0e-2 | 1.0e-4 | 2.6e-4 | 1.4e-3 |
|  | BKM (57+15) | BKM (88+15) | BPM (49) | BPM (65) |
| $\gamma = 2\sqrt{2}$ | 3.0e-2 | 8.0e-3 | 5.5e-4 | 3.2e-3 |

Table 2. $L_2$ norm of relative errors for 2D inhomogeneous convection-diffusion problems by the BKM and BPM

|  | BKM (33+15) | BKM (41+15) | BPM (25) | BPM (41) |
|---|---|---|---|---|
| $P^*=36$ | 8.4e-3 | 2.1e-4 | 9.0e-3 | 3.6e-4 |
|  | BKM (17+15) | BKM (25+15) | BPM (25) | BPM (41) |
| $P=540$ | 1.1e-3 | 3.5e-2 | 4.3e-3 | 4.1e-3 |

$^*P$ denotes Peclect number.

Table 3. $L_2$ norm of relative errors for 3D homogeneous Helmholtz problems by the BKM.

| Helmholtz ($\gamma = \sqrt{3}$) | | Helmholtz ($\gamma = 2\sqrt{3}$) | |
|---|---|---|---|
| 1.3e-2 (366) | 2.8e-3 (498) | 5.7e-2 (804) | 3.1e-3 (996) |

Table 4. $L_2$ norm of relative errors for 3D homogeneous convection-diffusion problems by the BKM.

| Convection-diffusion ($P=56$) | | Convection-diffusion ($P=560$) | |
|---|---|---|---|
| 7.0e-6 (114) | 3.5e-6 (174) | 2.2e-33 (114) | 2.4e-33 (174) |

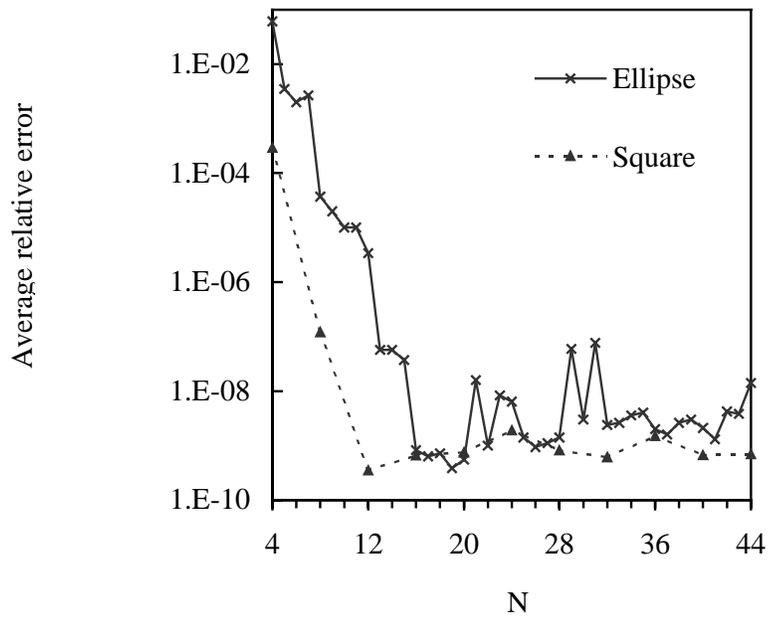

Fig. 3. Average relative error curves for Dirichlet convection-diffusion problem with 2D elliptical and square domains

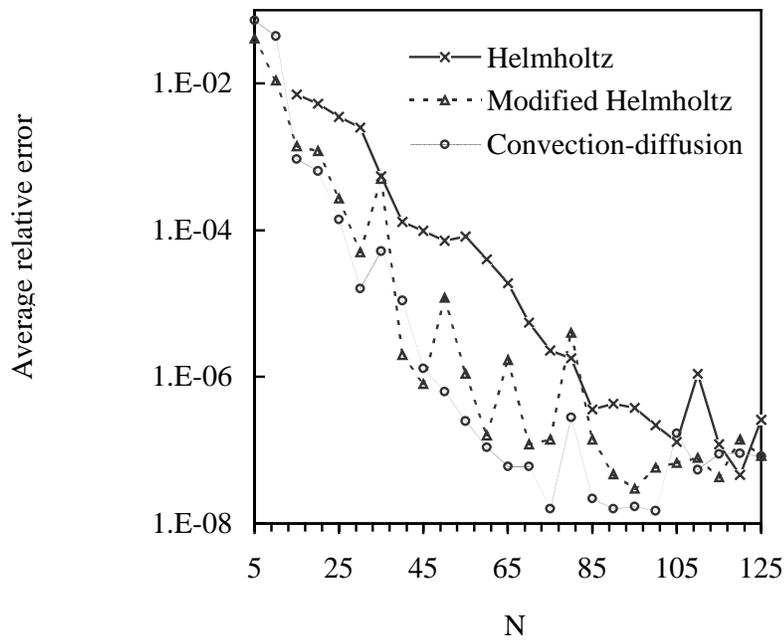

Fig. 4. Average relative error curve for Homogeneous Dirichlet Helmholtz, modified Helmholtz and convection-diffusion problems with 3D sphere domain

# 5. Modified Kansa method and its numerical validations

The Green integral solution of the previous PDE case is given by

$$u(x) = \int_\Gamma \left[ u \frac{\partial u^*(x,z)}{\partial n} - \frac{\partial u}{\partial n} u^*(x,z) \right] d\Gamma(z) + \int_\Omega f(z) u^*(x,z) d\Omega(z).$$

With a numerical integral scheme, we have

$$u(x) \cong \sum_{k=1}^{N+L} p(x, x_k) \frac{\partial u^*}{\partial n} u + \sum_{k=1}^{N+L} h(x, x_k) u^* \frac{\partial u}{\partial n} + \sum_{k=1}^{N+L} \omega(x, x_k) u^* f(x_k).$$

By analogy with the Fasshauer's Hermite scheme, we can construct

$$u(x) = \sum_{k=1}^{L_d} \alpha_k \varphi(r_k) + \sum_{k=L_d+1}^{L} \alpha_k \left[ -\frac{\partial \varphi(r_k)}{\partial n} \right] + \sum_{k=N+L+1}^{N+2L} \alpha_k \Re^* \{\varphi(r_k)\}.$$

Substituting the above expression into boundary and governing equations, we have the standard *Ax=b* formulation, where

$$A = \begin{bmatrix} \varphi & -\dfrac{\partial \varphi}{\partial n} & \Re^*\{\varphi\} \\ \dfrac{\partial \varphi}{\partial n} & -\dfrac{\partial^2 \varphi}{\partial n^2} & \dfrac{\partial \Re^*\{\varphi\}}{\partial n} \\ \Re\{\varphi\} & -\dfrac{\partial \Re\{\varphi\}}{\partial n} & \Re\Re^*\{\varphi\} \end{bmatrix}.$$

The above scheme is called the modified Kansa method (MKM) in contrast to the traditional Kansa method.

Table 5. $L_2$ norms of relative errors for Dirichlet Laplace and Helmholtz problems with a unit square domain by the MKM.

| Laplace | | Helmholtz ($\gamma = 2\sqrt{2}$) | |
|---|---|---|---|
| 8.7e-3 (49) | 1.4e-3 (81) | 1.5e-4 (25) | 1.5e-4 (36) |

The exact solution of Laplace problem is $u = 2\sin(\pi x) + \sin(2\pi y)$.

# 6. Remarks: merits and demerits

## Merits:

1. Very easy to learn and program.

2. Independent of geometric complexity and dimensionality, applicable to high-dimensional moving boundary problems.

3. Symmetric, meshfree, integration-free and spectral-convergence.

## Demerits:

1. Severe ill-conditioning of large dense RBF interpolation matrix

2. Immature mathematical theory: convergence, stability, and applicability.

3. Lacking rapid solution of global RBF interpolation of PDE's: localization and decomposition with preconditioning.